\documentclass[11pt]{article}

\usepackage{amssymb, amsmath, amsthm, graphicx}
\usepackage[left=1in,top=1in,right=1in]{geometry}

\date{}

\theoremstyle{plain}
      \newtheorem{theorem}{Theorem}[section]
      \newtheorem{lemma}[theorem]{Lemma}
            
      \newtheorem{problem}[theorem]{Problem}
      \newtheorem{observation}[theorem]{Observation}

      \newtheorem{conjecture}[theorem]{Conjecture}
\theoremstyle{definition}

\theoremstyle{remark}

	\newcommand{\RR}{{\mathbb R}}

\title{Bounded VC-dimension implies the Schur-Erd\H os conjecture}
\author{Jacob Fox\thanks{Stanford University, Stanford, CA. Supported by a Packard Fellowship and by NSF award DMS-1855635. Email: {\tt jacobfox@stanford.edu.}} \and J\'anos Pach\thanks{R\'enyi Institute of Mathematics, Budapest. Supported
by Swiss National Science Foundation Grants 200020-162884 and 200021-175977. Email:
{\tt pach@cims.nyu.edu}.}\and  Andrew Suk\thanks{Department of Mathematics, University of California at San Diego, La Jolla, CA, 92093 USA. Supported by NSF CAREER award DMS-1800746 and an Alfred Sloan Fellowship. Email: {\tt asuk@ucsd.edu}.} }

\begin{document}

\maketitle

\begin{abstract}
In 1916, Schur introduced the Ramsey number $r(3;m)$, which is the minimum integer $n$ such that for any $m$-coloring of the edges of the complete graph $K_n$, there is a monochromatic copy of $K_3$.  He showed that $r(3;m) \leq O(m!)$, and a simple construction demonstrates that $r(3;m) \geq 2^{\Omega(m)}$.  An old conjecture of Erd\H os states that $r(3;m)  = 2^{\Theta(m)}$.  In this note, we prove the conjecture for $m$-colorings with bounded VC-dimension, that is, for $m$-colorings with the property that the set system $\mathcal{F}$ induced by the neighborhoods of the vertices with respect to each color class has bounded VC-dimension.

\end{abstract}

\section{Introduction}

Given $n$ points and $n$ lines in the plane, their incidence graph is a bipartite graph $G$ that contains no $K_{2,2}$ as a subgraph. By a theorem of Erd\H os~\cite{Er38} and K\H ov\'ari-S\'os-Tur\'an~\cite{KST54}, this implies that the number of incidences between the points and the lines is $O(n^{3/2})$. However, a celebrated theorem of Szemer\'edi and Trotter~\cite{SzT83} states that the actual number of incidences is much smaller, only $O(n^{4/3})$, and this bound is tight. There are many similar examples, where extremal graph theory is applicable, but does not yield optimal results. What is behind this curious phenomenon? In the above and in many other examples, the vertices of $G$ are, or can be associated with, points in a Euclidean space, and the fact whether two vertices are connected by an edge can be determined by evaluating a bounded number of polynomials in the coordinates of the corresponding points. In other words, $G$ is a {\em semi-algebraic graph} of bounded complexity. As was proved by the authors, in collaboration with Sheffer and Zahl~\cite{FPSSZ17}, for semi-algebraic graphs, one can explore the {\em geometric} properties of the polynomial surfaces, including separator theorems and the so-called polynomial method~\cite{Gu16}, to obtain much stronger results for the extremal graph-theoretic problems in question. In particular, every $K_{2,2}$-free semi-algebraic graph of $n$ vertices with the complexity parameters associated with the point-line incidence problem has $O(n^{4/3})$ edges. This implies the Szemer\'edi-Trotter theorem.
\smallskip

There is a fast growing body of literature demonstrating that many important results in extremal combinatorics can be substantially improved, and several interesting conjectures proved, if we restrict our attention to semi-algebraic graphs and hypergraphs; see, e.g.,~\cite{APPRS05, FGLNP12, FPS16}. It is a major unsolved problem to decide whether this partly algebraic and partly geometric assumption can be relaxed and replaced by a purely combinatorial condition. A natural candidate is that the graph has {\em bounded Vapnik-Chervonenkis dimension} (in short, VC-dimension). The {\em VC-dimension of a set system} (hypergraph) $\mathcal{F}$ on the ground set $V$ is the {\em largest} integer $d$ for which there exists a $d$-element set $S\subset V$ such that for every subset $B\subset S$, one can find a member $A\in \mathcal{F}$ with $A\cap S=B$. The {\em VC-dimension of a graph} $G=(V,E)$ is the VC-dimension of the set system formed by the neighborhoods of the vertices, where the neighborhood of $v\in V$ is $N(v)=\{u\in v : uv\in E\}$. The VC-dimension, introduced by Vapnik and Chervonenkis~\cite{VC}, is one of the most useful combinatorial parameters that measures the complexity of graphs and hypergraphs. It proved to be relevant in many branches of pure and applied mathematics, including statistics, logic, learning theory, and real algebraic geometry. It has completely transformed combinatorial and computational geometry after its introduction to the subject by Haussler and Welzl~\cite{HaW87} in 1987. But can it be applied to extremal graph theory problems?
\smallskip

At first glance, this looks rather unlikely. Returning to our initial example, it is easy to verify that the Vapnik-Chervonenkis dimension of every $K_{2,2}$-free graph is at most $2$. Therefore, we cannot possibly improve the $O(n^{3/2})$ upper bound on the number of edges of $K_{2,2}$-free graphs by restricting our attention to graphs of bounded VC-dimension. Yet the goal of the present note is to solve the {\em Schur-Erd\H os problem}, one of the oldest open questions in Ramsey theory, by placing this restriction.
\smallskip

To describe the problem, we need some notation. For integers $k \geq 3$ and $m \geq 2$, the \emph{Ramsey number} $r(k;m)$ is the smallest integer $n$ such that any $m$-coloring of the edges of the complete $n$-vertex graph contains a monochromatic copy of $K_k$.  For the special case when $k = 3$, Issai Schur \cite{S} showed that
        $$\Omega(2^m) \le r(3;m) \le O(m!).$$
While the upper bound has remained unchanged over the last 100 years, the lower bound was successively improved. The current record is due to Xiaodong \emph{et al.}~\cite{fred} who showed that $r(3;m) \geq \Omega(3.199^m)$.  It is a major open problem in Ramsey theory to close the gap between the lower and upper bounds for $r(3;m)$.  Erd\H os \cite{erdos} offered cash prizes for solutions to the following problems.

\begin{conjecture}[\$100]\label{conjecture}
$$\lim\limits_{m \rightarrow \infty} (r(3;m))^{1/m}<\infty$$
\end{conjecture}

It was shown by Chung \cite{C} that $r(3;m)$ is supermultiplicative, so that the above limit exists.

\begin{problem}[\$250]
\hskip 0.4cm Determine \hskip 0.7cm $\lim\limits_{m \rightarrow \infty} (r(3;m))^{1/m}$
\end{problem}

It will be more convenient to work with the {\em dual VC-dimension}.   The \emph{dual} of a set system $\mathcal{F}$ is the set system $\mathcal{F}^{\ast}$ obtained by interchanging the roles of $V$ and $\mathcal{F}$. That is, the ground set of $\mathcal{F}^{\ast}$ is $\mathcal{F}$, and
$$\mathcal{F}^{\ast} = \{\{A\in \mathcal{F}: v\in A\}:v \in V\}.$$
We say that $\mathcal{F}$ has \emph{dual VC-dimension} $d$ if $\mathcal{F}^{\ast}$ has VC-dimension $d$.  Notice that $(\mathcal{F}^{\ast})^{\ast} = \mathcal{F}$, and it is known that if $\mathcal{F}$ has VC-dimension $d$, then $\mathcal{F}^{\ast}$ has VC-dimension at most $2^{d+1} - 1$ (see \cite{mat}). In particular, the VC-dimension of $\mathcal{F}$ is bounded if and only if the dual VC-dimension is.
\smallskip

Let $\chi$ be an $m$-coloring of the edges of the complete graph $K_n$ with colors $q_1,\ldots, q_m$, and let $V$ be the vertex set of $K_n$.  For $v \in V$ and $i \in [m]$, let $N_{q_i}(v) \subset V$ denote the neighborhood of $v$ with respect to the edges colored with color $q_i$.  We say that $\chi$ has \emph{VC-dimension} (or {\em dual VC-dimension}) $d$ if the set system $\mathcal{F} = \{N_{q_i}(v): i \in [m], v\in V\}$ has VC-dimension (resp., dual VC-dimension) $d$.

For $ k\geq 3$, $m\geq 2$, and $d \geq 2$, let $r_d(k;m)$ be the smallest integer $n$ such that for any $m$-coloring $\chi$ of the edges of $K_n$ with dual VC-dimension at most $d$ contains a monochromatic clique of size $k$.  Even for $m$-colorings with dual VC-dimension 2, we have $r_2(3;m) = 2^{\Omega(m)}$. Indeed, recursively take two disjoint copies of $K_{2^{m-1}}$, each of which is $(m-1)$-colored with dual VC-dimension at most $2$ and no monochromatic copy of $K_3$. Color all edges between these complete graphs with the $m$th color, to obtain an $m$-colored complete graph $K_{2^m}$ with the desired properties. Our main result shows that, apart from a constant factor in the exponent, this construction is tight.

\begin{theorem}\label{main}
For every $k \geq 3$ and $d \geq 2$, there is a constant $c = c(k,d)$ such that $r_d(k;m) \leq 2^{cm}.$  In other words, for every $m$-coloring of the edges of a complete graph of $2^{cm}$ vertices with dual VC-dimension $d$, there is a monochromatic complete subgraph of $k$ vertices.
\end{theorem}

It follows from the Milnor-Thom theorem~\cite{Mi64} (proved 15 years earlier by Ole\u inik and Petrovskii~\cite{PeO49}) that every $m$-coloring of the ${n\choose 2}$ pairs induced by $n$ points in $\RR^d$, which is {\em semi-algebraic} with bounded complexity, has bounded VC-dimension and, hence, bounded dual VC-dimension.

In a recent paper~\cite{fps}, we proved Conjecture~\ref{conjecture} for semi-algebraic $m$-colorings of bounded complexity. Our proof heavily relied on the topology of Euclidean spaces: it was based on the cutting lemma of Chazelle {\em et al.}~\cite{chaz} and vertical decomposition. These arguments break down in the combinatorial setting, for $m$-colorings of bounded VC-dimension. In what follows, instead of using ``regular'' space decompositions with respect to a set of polynomials, our main tool will be a partition result for abstract hypergraphs, which can be easily deduced from the dual of Haussler's packing lemma~\cite{H95}. The proof of this partition result will be given in Section~\ref{sec2}, while Section~\ref{pf} contains the proof of Theorem~\ref{main}.

\smallskip

To simplify the presentation, throughout this paper we omit the floor and ceiling signs whenever they are not crucial. All logarithms are in base 2.

\section{A partition lemma}\label{sec2}

 Let $\mathcal{F}$ be a set system with dual VC-dimension $d$ and with ground set $V$. Given two points $u,v \in V$, we say that a set $A\in \mathcal{F}$ \emph{crosses} the pair  $\{u,v\}$ if $A$ contains at least one member of $\{u,v\}$, but not both.  We say that the set $X\subset V$ is \emph{$\delta$-separated} if for any two points $u,v \in X$, there are at least $\delta$ sets in $\mathcal{F}$ that cross the pair $\{u,v\}$.   The following \emph{packing lemma} was proved by Haussler in \cite{H95}.

\begin{lemma}\label{dual}

Let $\mathcal{F}$ be a set system on a ground set $V$ such that $\mathcal{F}$ has dual VC-dimension $d$.   If $X\subset V$ is $\delta$-separated, then $|X| \leq c_1(|\mathcal{F}|/\delta)^{d}$ where $c_1= c_1(d)$.

\end{lemma}

As an application of Lemma \ref{dual}, we obtain the following partition lemma.

\begin{lemma}\label{partition}

Let $\mathcal{F}$ be a set system on a ground set $V$ with dual VC-dimension $d$.

Then there is a constant $c_2 = c_2(d)$ such that for any $\delta$ satisfying $1 \leq \delta
\leq |\mathcal{F}|$, there is a partition $V = S_1\cup \cdots \cup S_r$ of $V$ into
$r \leq c_2(|\mathcal{F}|/\delta)^{d}$ parts, each of size at most $\frac{2n}{c_1(|\mathcal{F}|/\delta)^{d}}$, such that any pair of vertices from the same part $S_t$ is crossed by at most $2\delta$ members of $\mathcal{F}$. (Here $c_1=c_1(d)$ is the same constant as in Lemma \ref{dual}.

\end{lemma}

\begin{proof}

Let $X = \{x_1,\ldots, x_{r'}\}$ be a maximal subset of $V$ that is $\delta$-separated with respect to $\mathcal{F}$.  By Lemma \ref{dual}, $|X| = r' \leq c_1(|\mathcal{F}|/\delta)^{d}$. We define a partition $V = S_1\cup \cdots \cup S_{r'}$ of the vertex set such that $v \in S_i$ if $i$ is the smallest index such that the number of sets from $\mathcal{F}$ that cross the pair $\{v,x_i\}$ is at most $\delta$.  Such an $i$ always exists since $X$ is maximal.   By the triangle inequality, for any two vertices $u,v \in S_i$, there are at most $2\delta$ sets in $\mathcal{F}$ that cross the pair $\{u,v\}$.

If a part $S_i$ has size more than $\frac{2n}{c(|\mathcal{F}|/\delta)^{d}}$, we partition $S_i$ (arbitrarily) into parts of size $\left\lfloor \frac{2n}{c_1(|\mathcal{F}|/\delta)^{d}}\right\rfloor$ and possibly one additional part of size less than $\left\lfloor\frac{2n}{c_1(|\mathcal{F}|/\delta)^{d}}\right\rfloor$.  Let $\mathcal{P}:V= S_1\cup \cdots \cup S_r$ be the resulting partition, where $r \leq c_2(|\mathcal{F}|/\delta)^{d}$ and $c_2 = c_2(d)$.  Then $\mathcal{P}$ satisfies the above properties.  \end{proof}

\section{Proof of Theorem \ref{main}}\label{pf}

 Let $d,k_1,\ldots,k_m$ be positive integers.  We define the Ramsey number $r_d(k_1,\ldots, k_m)$ to be the smallest integer $n$ with the following property. For any $m$-coloring $\chi$ of the edges of $K_n$ with colors $\{q_1,\ldots, q_m\}$ such that $\chi$ has dual VC-dimension at most $d$, there is a monochromatic copy of $K_{k_i}$ in color $q_i$ for some $1 \leq i \leq m$.  We now prove the following theorem, from which Theorem \ref{main} immediately follows.

\begin{theorem}\label{gen}

For fixed integers $d,k\geq 1$, if $k_1,\ldots, k_m \leq k$, then $r_d(k_1,\ldots, k_m) = 2^{O(m)}.$

\end{theorem}

\begin{proof}

Let $c = c(d,k)$ be a large constant that will be determined later.  We will show that $r_d(k_1,\ldots, k_m) \leq 2^{c\sum_{i = 1}^mk_i}$ by induction on $s = \sum_{i = 1}^mk_i$.  The base case $s \leq k2^{16dk}$ follows by setting $c$ to be sufficiently large.

For the inductive step, assume that $s > k2^{16dk}$ and that the statement holds for all $s' < s$.  Thus, we have $m\geq 2^{16dk}$.  Let $n = 2^{cs}$ and let $\chi$ be an $m$-coloring of the edges of $K_n$ with colors $q_1,\ldots, q_m$ such that the set system $\mathcal{F} = \{N_{q_i}(v): v \in V(K_n), i \in [m]\}$ has dual VC dimension at most $d$.

Set $\mathcal{F}_0 = \mathcal{F}$ and $V_0 = V$, and let $\log^{(j)} m$ denote the $j$-fold iterated logarithm function, where $\log^{(0)} m = m$ and $\log^{(j)} m = \log(\log^{(j - 1)} m )$.  For $j\geq 1$ such that $\log^{(j)}m > 2^{8dk}$, we will recursively construct a set system $\mathcal{F}_j$, whose ground set is $V_j \subset V$, such that

\begin{enumerate}

\item $\mathcal{F}_j = \{N_{q_i}(v)\cap V_j: v \in V_j, q_i \in Q_v\},$ where $Q_v\subset \{q_1,\ldots, q_m\}$ and $|Q_v| \leq \log^{(j)}m$.  Hence, $|\mathcal{F}_j| \leq n\log^{(j)}m$.

\item $|V_j| \geq n - \frac{n}{\log^{(j-1)} m}$.

\item $\mathcal{F}_j$ covers at least ${n\choose 2} - \frac{8n^2}{\log^{(j-1)}m}$ edges of $K_n$, where an edge $uv \in E(K_n)$ is \emph{covered} by $\mathcal{F}_j$ if $(N_{q_i}(u)\cap V_j),(N_{q_i}(v)\cap V_j) \in \mathcal{F}_j$ where $q_i = \chi(uv)$.

\end{enumerate}

Having obtained $\mathcal{F}_{j}$ and $V_{j}$ with the properties described above, we obtain $\mathcal{F}_{j+1}$ and $V_{j+1}$ as follows.  Let $B_j\subset E(K_n)$ denote the set of edges that are not covered by $\mathcal{F}_j$.  Hence, $|B_j| \leq \frac{8n^2}{\log^{(j-1)}m}$.  We apply Lemma \ref{partition} to $\mathcal{F}_{j}$, whose ground set is $V_j$, with parameter $\delta = \frac{|\mathcal{F}_{j}|}{(\log^{(j)} m)^4}$, and obtain a partition $\mathcal{P}:V_{j}  = S_{1}\cup \cdots \cup S_{r}$, where $r  \leq c_2(\log^{(j)} m)^{4d}$ and $c_2$ is defined in Lemma \ref{partition}, such that $\mathcal{P}$ has the properties described in Lemma \ref{partition}.  For each part $S_{t} \in \mathcal{P}$, let $Q_{t} \subset \{q_1,\ldots, q_m\}$ be the set of colors such that $q_i \in Q_{t}$ if there is a vertex $v\in S_{t}$ such that

$$|\{u \in V_j: \chi(uv) = q_i, uv \not\in B_j\}|\geq \frac{n}{(\log^{(j)} m)^2}.$$

\noindent  Let $Q'_{t}  \subset \{q_1,\ldots, q_m\}$ be the set of colors such that $q_i \in Q'_{t}$ if the vertex set $S_{t}$ contains a monochromatic copy of $K_{k_i - 1}$ in color $q_i$.

\begin{observation}

If there is a color $q_i \in Q_{t}\cap Q'_{t}$, then $\chi$ produces a monochromatic copy of $K_{k_i}$ in color $q_i$.

\end{observation}

\begin{proof}

Suppose $q_i \in Q_{t}\cap Q'_{t}$ and let $X  = \{x_1,\ldots, x_{k_i - 1}\} \subset S_{t}$ be the vertex set of a monochromatic clique of order $k_i - 1$ in color $q_i$.  Fix $v \in S_{t}$ such that for $U  = \{u \in V_j: \chi(uv) = q_i, uv \not\in B_j\}$, we have $|U|\geq \frac{n}{(\log^{(j)} m)^2}.$   Notice that if $X\not \subset (N_{q_i}(u)\cap V_j)$, where $u\in U$, then the set $(N_{q_i}(u)\cap V_j)$ crosses the pair $\{x,v\}$ for some $x \in X$.  Moreover,  $(N_{q_i}(u)\cap V_j) \in \mathcal{F}_{j}$ since $uv \not\in B_j$. Since there are at most $2\delta = \frac{2|\mathcal{F}_{j}|}{(\log^{(j)} m)^4}$ sets in $\mathcal{F}_{j}$ that cross $\{x,v\}$, there are at most $\frac{2k|\mathcal{F}_{j}|}{(\log^{(j )} m)^4}$ sets in $\{ N_{q_i}(u)\cap V_j: u \in U\}\subset \mathcal{F}_j$ that do not contain $X$. On the other hand, $$|U| - k_i \geq \frac{n}{(\log^{(j)} m)^2} - k_i > \frac{2k|\mathcal{F}_{j}|}{(\log^{(j )} m)^4},$$ where the last inequality follows from the fact that $|\mathcal{F}_{j}| \leq n\log^{(j)}m$ and $\log^{(j)}m >2^{8dk}$.  Hence, there must be a neighborhood $(N_{q_i}(u)\cap V_j)$ that contains $X$, which implies that $X\cup \{u\}$ induces a monochromatic copy of $K_{k_i}$ in color $q_i$. \end{proof}

By the observation above, we can assume that $Q_{t}\cap Q'_{t} = \emptyset$ for every $t$, since otherwise we would be done.

\begin{observation}

If there is a part $S_{t} \in \mathcal{P}$ such that $|S_{t}| \geq n/(\log^{(j)} m)^{6d}$ and $|Q_{t}| \geq \log^{(j + 1)} m$, then $S_t$ contains a monochromatic copy of $K_{k_i}$ in color $q_i$ where $q_i\in Q'_t$.

\end{observation}

\begin{proof}  For sake of contradiction, suppose $S_t \in \mathcal{P}$ does not contain a monochromatic copy of $K_{k_i}$ in color $q_i \in Q'_t$.  Since $Q_{t}\cap Q'_{t} = \emptyset$, $S_t$ also does not contain a monochromatic copy of $K_{k_i -1}$ in color $q_i \in Q_t$.  So if $|Q_{t}| \geq \log^{(j + 1)} m$, we have $|Q'_t| \leq m - \log^{(j + 1)} m$.  By the induction hypothesis, we have

$$\frac{n}{(\log^{(j )} m)^{6d}} \leq |S_{t}| < 2^{c(s - \log^{(j  + 1)} m)}.$$

\noindent Since $c = c(d,k)$ is sufficiently large, we have $n  < 2^{cs}$ which is a contradiction.  \end{proof}

\noindent Hence, we can assume that for each part $S_{t}\in \mathcal{P}$ such that $|S_{t}| \geq n/(\log^{(j)} m)^{6d}$, we have $|Q_{t}| < \log^{(j  + 1)} m$.

We now define $\mathcal{F}_{j+1}$ and $V_{j+1}$ as follows.  Start with the set system $\mathcal{F}_{j}$ whose ground set is $V_{j}$. For each vertex $v \in V_{j}$ that lies in a part $S_t\in \mathcal{P}$ with $|S_{t}| <n/(\log^{(j )} m)^{6d}$, we remove all sets in $\mathcal{F}_{j}$ of the form $N_{q_i}(v)\cap V_j$, that is, we remove all neighborhoods generated by $v$ in $\mathcal{F}_j$.   For each vertex $v \in S_{t}$, such that $S_t \in \mathcal{P}$ and $|S_t| \geq n/(\log^{(j )} m)^{6d}$, we remove all sets of the form $N_{q_i}(v) \cap V_j \in \mathcal{F}_{j}$ if $q_i \not\in Q_{t}$.  Let $\mathcal{F}_{j+1}$ be the remaining set system induced on the ground set $V_{j+1}$, where $V_{j+1}\subset V_j$ is the set of vertices obtained by deleting all parts $S_t\in \mathcal{P}$ such that $|S_{t}| <n/(\log^{(j)} m)^{6d}$.   Clearly, each vertex $v \in V_{j+1}$ only contributes at most $|Q_{t}| \leq \log^{(j+1)} m$ sets in $\mathcal{F}_{j+1}$, and each vertex $v \in V\setminus V_{j+1}$ does not contribute any sets in $\mathcal{F}_{j + 1}$.   Hence, $|\mathcal{F}_{j+1}| \leq n \log^{(j+1)} m$.  Moreover,

$$\begin{array}{ccl}
   |V_{j+1}| & \geq  & |V_{j}| -  c_2(\log^{(j)} m)^{4d}\frac{n}{(\log^{(j )} m)^{6d}} \\\\
     & \geq  & n - \frac{n}{\log^{(j-1)}m} -  \frac{c_2n}{(\log^{(j )} m)^{2d}}\\\\
     & \geq  &  n - \frac{n}{\log^{(j)} m}.
 \end{array}$$

Finally, it remains to show that $\mathcal{F}_{j+1}$ covers at least ${n\choose 2} - \frac{8n^2}{\log^{(j)} m}$ edges of $K_n$.  Let $B_{j + 1}\subset E(K_n)$ denote the set of edges that are not covered by $\mathcal{F}_{j + 1}$.  If $uv\in B_{j + 1}$, then either

\begin{enumerate}

\item $uv \in B_{j}$, or

\item $u$ (or $v$) lies inside a part $S_t \in \mathcal{P}$ such that $|S_t| \leq\frac{n}{(\log^{(j)}m)^{6d}}$, or

\item both $u$ and $v$ lie inside the same part $S_t \in \mathcal{P}$, or

\item $uv$ is covered by $\mathcal{F}_j$, but is not covered by $\mathcal{F}_{j + 1}$ since $v \in S_t \in \mathcal{P}$ and $\chi(u,v) \not\in Q_{t}$.

\end{enumerate}

\noindent By assumption,

\begin{equation}\label{one}
|B_j| \leq \frac{8n^2}{\log^{(j -1)} m}.
\end{equation}

\noindent The number of edges of the second type is at most

\begin{equation}\label{two}
\frac{n^2}{(\log^{(j)}m)^{6d}}.
\end{equation}

\noindent The number of edges of the third type is at most

\begin{equation}\label{three}
\sum\limits_{i = 1}^r{|S_t| \choose 2} \leq c_2(\log^{(j)} m)^{4d} \left(\frac{2n}{c_1(\log^{(j)}m)^{4d}}\right)^2 = \frac{4c_2n^2}{(c_1)^2(\log^{(j)}m)^{4d}},
\end{equation}

\noindent where $c_1$ is defined in Lemma \ref{dual}. Finally, let us bound the number of edges of the fourth type.  Fix $v \in S_t \in \mathcal{P}$ such that $|S_t| > \frac{n}{(\log^{(j)}m)^{6d}}$, and let us consider all edges incident to $v$ that are covered by $\mathcal{F}_j$.  Since  $v$ contributed at most $\log^{(j)}m$ sets in $\mathcal{F}_j$, there are at most $\log^{(j)}m$ distinct colors among these edges.  Fix such a color $q_i$ such that $q_i \not\in Q_t$, and consider the set of vertices

$$U = \{u \in V_{j + 1}: \chi(uv) = q_i,uv\not\in B_j\}.$$

\noindent By definition of $Q_t$, we have $|U| < \frac{n}{(\log^{(j)} m)^2}$.   Therefore, the number of edges incident to $v$ of the fourth type is at most $$\log^{(j)}m \frac{n}{(\log^{(j)} m)^2} =  \frac{n}{\log^{(j)} m}.$$  Hence, the total number of edges of the fourth type is at most

\begin{equation}\label{four}
 \frac{n^2}{\log^{(j)} m}.
\end{equation}

\noindent Thus by summing (\ref{one}), (\ref{two}), (\ref{three}), (\ref{four}), and using the fact that $\log^{(j)}m > 2^{8dk}$, we have

$$|B_{j + 1}| \leq \frac{8n^2}{\log^{(j-1)} m} + \frac{n^2}{(\log^{(j)}m)^{6d}} + \frac{4c_2n^2}{(c_1)^2(\log^{(j)} m)^{4d}} + \frac{n^2}{ \log^{(j)} m}< \frac{8n^2}{\log^{(j)}m }.$$

\noindent Hence, $\mathcal{F}_{j + 1}$ covers at least ${n\choose 2} - \frac{8n^2}{\log^{(j)}m }$ edges of $K_n$.

Let $w$ be the minimum integer such that $\log^{(w)} m < 2^{8dk}$.  Then we have $\mathcal{F}_{w}, V_{w}, B_{w}$ with the properties described above, where $B_w \subset E(K_n)$ is the set of edges not covered by $\mathcal{F}_w$.  This implies $|B_{w}| \leq n^2/2^{8dk} < n^2/8$ and $|V_{w}| \geq 7n/8$.  Since

$${7n/8 \choose 2} - \frac{n^2}{8} \geq \frac{n^2}{4},$$

\noindent an averaging argument shows that there is a vertex $v \in V_{w}$ that is incident to at least $n/2$ edges that are covered by $\mathcal{F}_{w}$.  Since $v$ contributes at most $\log^{(w)}m < 2^{8dk}$ sets in $\mathcal{F}_{w}$, there is a color $q_i$ such that $$|N_{q_i}(v)| \geq \frac{n}{2\cdot 2^{8dk}} \geq 2^{c(s-1)},$$ where the second inequality follows from the fact that $c = c(d,k)$ is sufficiently large. Therefore, by induction, the set $N_{q_i}(v) \subset V$ contains a monochromatic copy of $K_{k_{\ell}}$ in color $q_{\ell} \in \{q_1,\ldots, q_m\}\setminus q_i$, in which case we are done, or contains a monochromatic copy of $K_{k_i - 1}$ in color $q_i$.  In the latter case, we obtain a monochromatic $K_{k_i}$ in color $q_i$ by including vertex $v$.  This completes the proof of Theorem \ref{gen}.\end{proof}

\section{Concluding remarks}

We have established tight bounds for multicolor Ramsey numbers for graphs with bounded VC-dimension.  It would be interesting to prove other well-known conjectures in extremal graph theory for graphs and hypergraphs with bounded VC-dimension, especially the notorious Erd\H os-Hajnal conjecture.

An old result of Erd\H os and Hajnal \cite{EH} states that for every hereditary property $P$ which is not satisfied by all graphs,
there exists a constant $\varepsilon(P) > 0$ such that every graph of $n$ vertices with property $P$ has a clique or
an independent set of size at least $e^{\varepsilon(P)\sqrt{\log n}}$.  They conjectured that this bound can be improved to $n^{\varepsilon(P)}$.  Thus, every graph $G$ on $n$ vertices with bounded VC-dimension contains a clique or an independent set of size $e^{\Omega(\sqrt{\log n})}$.  In \cite{FPSEH}, the authors improved this bound to $e^{(\log n)^{1 - o(1)}}$.  However, the following conjecture remains open.

\begin{conjecture}
For $d\geq 2$, there exists a constant $\varepsilon(d)$ such that every graph on $n$ vertices with VC-dimension at most $d$ contains a clique or an independent set of size $n^{\varepsilon(d)}$.

\end{conjecture}

\end{document}